\magnification=1200
%%%%%%%%%%%%%%%%%%%%%%%%%%%%%%%%%%%%%%%%%%%%%%%%%%%%%%%%%%%%%%%%%
%                       \input zaj.sty                          %
\input amssym.def
\input amssym.tex

\hsize=36truecc

\font\secbf=cmb10 scaled 1200
\font\eightrm=cmr8
\font\sixrm=cmr6

\font\eighti=cmmi8

\font\sixi=cmmi6
\skewchar\eighti='177 \skewchar\sixi='177

\font\eightsy=cmsy8
\font\sixsy=cmsy6
\skewchar\eightsy='60 \skewchar\sixsy='60

\font\eightit=cmti8

\font\eightbf=cmbx8
\font\sixbf=cmbx6

\let\sc=\tensc

\font\eightsc=cmcsc10 scaled 800

             % roman do streszcze¤ i literatury
              % italic do literatury

    % Font na tytu'y
 % Font na podtytu'y
\font\secbf=cmb10 scaled 1200
\font\subsecfont=cmb10 scaled \magstephalf
%%%%%%%%%%%%%%%%%%%%%%%%%%%%%%%%%%%%%%%%%%%%%%%%%%
\font\amb=cmmib10

\font\ambi=cmmib10 scaled 700

\newfam\mbfam \def\mb{\textfont1=\amb\fam\mbfam\amb\scriptfont1=\ambi}

\textfont\mbfam\amb \scriptfont\mbfam\ambi

\def\bm#1{\mathchoice
{\hbox{\mb\textfont1=\amb$#1$}}%
{\hbox{\mb\textfont1=\amb$#1$}}%
{\hbox{\mb$\scriptstyle\textfont1=\ambi#1$}}%
{\hbox{\mb$\scriptscriptstyle\textfont1=\ambi#1$}}}
%%%%%%%%%%%%%%%%%%%%%%%%%%%%%%%%%%%%%%%%%%%%%%%%%%%

\def\aa{\def\rm{\fam0\eightrm}%
  \textfont0=\eightrm \scriptfont0=\sixrm \scriptscriptfont0=\fiverm
  \textfont1=\eighti \scriptfont1=\sixi \scriptscriptfont1=\fivei
  \textfont2=\eightsy \scriptfont2=\sixsy \scriptscriptfont2=\fivesy
  \textfont3=\tenex \scriptfont3=\tenex \scriptscriptfont3=\tenex
  \def\sc{\eightsc}
  \def\it{\fam\itfam\eightit}%
  \textfont\itfam=\eightit
  \def\bf{\fam\bffam\eightbf}%
  \textfont\bffam=\eightbf \scriptfont\bffam=\sixbf
   \scriptscriptfont\bffam=\fivebf
  \normalbaselineskip=9.7pt
  \setbox\strutbox=\hbox{\vrule height7pt depth2.6pt width0pt}%
  \normalbaselines\rm}

\def\Proof{\vskip12pt\noindent{\bf Proof.} }

\def\Remark#1{\vskip12pt\noindent{\bf Remark #1}}

\def\m@th{\mathsurround=0pt}

\def\cc#1{\hbox to .89\hsize{$\displaystyle\hfil{#1}\hfil$}\cr}
\def\lc#1{\hbox to .89\hsize{$\displaystyle{#1}\hfill$}\cr}
\def\rc#1{\hbox to .89\hsize{$\displaystyle\hfill{#1}$}\cr}

\def\eqal#1{\null\,\vcenter{\openup\jot\m@th
  \ialign{\strut\hfil$\displaystyle{##}$&&$\displaystyle{{}##}$\hfil
      \crcr#1\crcr}}\,}

\def\section#1{\vskip 22pt plus6pt minus2pt\penalty-400
        {{\secbf
        \noindent#1\rightskip=0pt plus 1fill\par}}
        \par\vskip 12pt plus5pt minus 2pt
        \penalty 1000}

\def\subsection#1{\vskip 20pt plus6pt minus2pt\penalty-400
        {{\subsecfont
        \noindent#1\rightskip=0pt plus 1fill\par}}
        \par\vskip 8pt plus5pt minus 2pt
        \penalty 1000}

\def\subsubsection#1{\vskip 18pt plus6pt minus2pt\penalty-400
        {{\subsecfont
        \noindent#1}}
        \par\vskip 7pt plus5pt minus 2pt
        \penalty 1000}

\def\center#1{{\begingroup \leftskip=0pt plus 1fil\rightskip=\leftskip
\parfillskip=0pt \spaceskip=.3333em \xspaceskip=.5em \pretolerance 9999
\tolerance 9999 \parindent 0pt \hyphenpenalty 9999 \exhyphenpenalty 9999
\par #1\par\endgroup}}

\def\\{\hfill\break}

\def\kwadrat{\hfill$\square$}
\def\mida#1{{{\null\kern-4.2pt\left\bracevert\vbox to 6pt{}\!\hbox{$#1$}\!\right\bracevert\!\!}}}
\def\midy#1{{{\null\kern-4.2pt\left\bracevert\!\!\hbox{$\scriptstyle{#1}$}\!\!\right\bracevert\!\!}}}

\def\diagint{{\raise1.5pt\hbox{$\scriptscriptstyle\diagup$}\hskip-8.7pt\intop}}

\def\divv{{\rm div}\,}
\def\rot{{\rm rot}\,}
\def\const{{\rm const}}

\def\today{${\scriptscriptstyle\number\day-\number\month-\number\year}$}
\footline={{\hfil\rm\the\pageno\hfil${\scriptscriptstyle\rm\jobname}$\ \ \today}}
%%%%%%%%%%%%%%%%%%%%%%%%%%%%%%%%%%%%%%%%%%%%%%%%%%%%%%%%%%%%%%%%%%%%%%%%%%%%%%

\hsize=36truecc
\baselineskip=16truept

\def\D{{\bm{D}}}

\def\I{{\bm{I}}}

\def\T{{\bm{T}}}

\def\0{{\bf0}}

\def\sup{\mathop{\rm sup}\limits}

\def\\{\hfil\break}
\def\N{{\Bbb N}}
\def\R{{\Bbb R}}

\center{\secbf Long time estimate of solutions to 3d Navier-Stokes 
equations coupled with the heat convection}

\vskip3cm
\centerline{\bf Jolanta Soca\l a$^1$, Wojciech M. Zaj\c aczkowski$^2$}
\ 

\item{$^1$} State Higher Vocational School in Racib\'orz,\\
S\l owacki Str. 55, 47-400 Racib\'orz, Poland,\\
e-mail: jolanta\_socala@interia.pl
\item{$^2$}Institute of Mathematics, Polish Academy of Sciences\\
\'Sniadeckich 8, 00-956 Warsaw, Poland\\
and Institute of Mathematics and Cryptology, Cybernetics Faculty,\\
Military University of Technology,\\
Kaliskiego 2, 00-908 Warsaw, Poland\\
e-mail:wz@impan.pl
\vskip1cm

\noindent
{\bf Abstract.} We examine the Navier-Stokes equations with homogeneous slip 
boundary conditions coupled with the heat equation with homogeneous Neumann 
conditions in a bounded domain in $\R^3$. The considered domain is a cylinder 
with $x_3$-axis. The aim of this paper is to show long time estimates without 
smallness of the initial velocity, the initial temperature and the external 
force. To prove the estimate we need however smallness of $L_2$ norms 
of derivatives with respect to $x_3$ of the initial velocity, the initial 
temperature and the external force.



\vfil\eject

\section{1. Introduction}

The aim of this paper is to derive long time a priori estimate for some 
initial-boundary value problem for a system of the Navier-Stokes equations 
coupled with the heat equation. We assume the slip boundary conditions for the 
Navier-Stokes equations and the Neumann condition for the heat equations. 
We examine the problem in a straight finite cylinder. To obtain the estimate 
we follow the ideas from [5, 7, 8] and the considered solution remains close 
to a two-dimensional solution. The estimate is the first and the most 
important step to prove the existence of solutions to the problem (see (1.1)) 
by the Leray-Schauder fixed point theorem (see the next paper of the authors 
[6]).

\noindent
We consider the following problem
$$\eqal{
&v_{,t}+v\cdot\nabla v-\divv\T(v,p)=\alpha(\theta)f\quad &{\rm in}\ \ 
\Omega^T=\Omega\times(0,T),\cr
&\divv v=0\quad &{\rm in}\ \ \Omega^T,\cr
&\theta_{,t}+v\cdot\nabla\theta-\varkappa\Delta\theta=0\quad &{\rm in}\ \ 
\Omega^T,\cr
&\bar n\cdot\D(v)\cdot\bar\tau_\alpha=0,\ \ \alpha=1,2\quad &{\rm on}\ \ 
S^T=S\times(0,T),\cr
&\bar n\cdot\bar v=0\quad &{\rm on}\ \ S^T,\cr
&\bar n\cdot\nabla\theta=0\quad &{\rm on}\ \ S^T,\cr
&v|_{t=0}=v(0),\ \ \theta|_{t=0}=\theta(0)\quad &{\rm in}\ \ \Omega,\cr}
\leqno(1.1)
$$
where by $x=(x_1,x_2,x_3)$ we denote the Cartesian coordinates, 
$\Omega\subset\R^3$ is a cylindrical type domain parallel to the axis $x_3$ 
with arbitrary cross section, $S=\partial\Omega$, 
$v=(v_1(x,t),v_2(x,t),v_3(x,t))\in\R^3$ is the velocity of the fluid motion, 
$p=p(x,t)\in\R^1$ the pressure, $\theta=\theta(x,t)\in\R_+$ the temperature, 
$f=(f_1(x,t),f_2(x,t),f_3(x,t))\in\R^3$ the external force field, $\bar n$ is 
the unit outward normal vector to the boundary $S$, $\bar\tau_\alpha$, 
$\alpha=1,2$, are tangent vectors to $S$ and the dot denotes the scalar 
product in $\R^3$. We define the stress tensor by
$$
\T(v,p)=\nu\D(v)-p\I,
$$
where $\nu$ is the constant viscosity coefficient, $\I$ is the unit matrix 
and $\D(v)$ is the dilatation tensor of the form
$$
\D(v)=\{v_{i,x_j}+v_{j,x_i}\}_{i,j=1,2,3}.
$$
Finally $\varkappa$ is a positive heat conductivity coefficient. 

\noindent
We assume that $S=S_1\cup S_2$, where $S_1$ is the part of the boundary which 
is parallel to the axis $x_3$ and $S_2$ is perpendicular to $x_3$. Hence
$$
S_1=\{x\in\R^3:\ \varphi_0(x_1,x_2)=c_*,\ -b<x_3<b\}
$$
and
$$
S_2=\{x\in\R^3:\ \varphi_0(x_1,x_2)<c_*,\ x_3\ {\rm is\ equal\ either\ to}\ 
-b\ {\rm or}\ b\},
$$
where $b,c_*$ are positive given numbers and $\varphi_0(x_1,x_2)$ describes 
a sufficiently smooth closed curve in the plane $x_3=\const$. We can assume 
$\bar\tau_1=(\tau_{11},\tau_{12},0)$, $\bar\tau_2=(0,0,1)$ and 
$\bar n=(\tau_{12},-\tau_{11},0)$ on $S_1$.
Assume that $\alpha\in C^2(\R)$ and $\Omega^T$ satisfies the weak $l$-horn 
condition, where $l=(2,2,1)$ (see [2, Ch. 2, Sect. 8]).

\noindent 
Moreover we assume that $\Omega^T$ is not axially symmetric. Now we formulate 
the main result of this paper. Let $g=f_{,x_3}$, $h=v_{,x_3}$, $q=p_{,x_3}$, 
$\vartheta=\theta_{,x_3}$, $\chi=(\rot v)_3$, $f=(\rot f)_3$. Assume that 
$\|\theta(0)\|_{L_\infty(\Omega)}<\infty$. Define
$$
a:[0,\infty)\to[0,\infty),\quad a(x)=\sup\{|\alpha(y)|+|\alpha'(y)|:|y|\le x\}
$$
and
$$
a(\theta(x))\le c_1,
\leqno(1.1')
$$
where $c_1=a(\|\theta(0)\|_{L_\infty(\Omega)})$. The inequality $(1.1')$ is 
justified in view of Lemma 2.3, Remark 2.4 and properties of function $a(x)$. 
Moreover assume that ${5\over3}<\sigma<\infty$, ${5\over3}<\varrho<\infty$, 
${5\over\varrho}-{5\over\sigma}<1$ and for $t\le T$
\vskip6pt

\item{1.} $c_1\|g\|_{L_2(0,t;L_{6/5}(\Omega))}+
c_1c_0\|f\|_{L_\infty(0,t;L_3(\Omega))}\\
+c_1\|F\|_{L_2(0,t;L_{6/5}(\Omega))}+
c_1\|f_3\|_{L_2(0,t;L_{4/3}(S_2))}\\
+\|h(0)\|_{L_2(\Omega)}+\|\vartheta(0)\|_{L_2(\Omega)}+
\|\chi(0)\|_{L_2(\Omega)}+\psi(c_0)\\
+c_0^2(c_1\|f\|_{L_2(0,t;L_{6/5}(\Omega))}+
\|v(0)\|_{L_2(\Omega)})\le k_1<\infty$,

\item{2.} $\|f\|_{L_2(0,t;L_3(\Omega))}\le k_2<\infty$,

\item{3.} $\|f\|_{L_2(\Omega^t)}+\|v(0)\|_{H^1(\Omega)}\le k_3<\infty$,

\item{4.} $c_1\|f\|_{L_\infty(\Omega^t)}e^{cc_1^2k_2^2}k_1+c_1
\|g\|_{L_\sigma(\Omega^t)}+\|\vartheta(0)\|_{W_\sigma^{2-2/\sigma}(\Omega)}\\
+\|h(0)\|_{W_\sigma^{2-2/\sigma}(\Omega)}\le k_4<\infty$,

\item{5.} $c_1\|g\|_{L_2(0,t;L_{6/5}(\Omega))}+c_1
\|f_3\|_{L_2(0,t;L_{4/3}(S_2))}\\
+\|h(0)\|_{L_2(\Omega)}+\|\vartheta(0)\|_{L_2(\Omega)}\le d<\infty$

\item{6.} $c_1+\|f\|_{L_\varrho(\Omega^t)}+
\|v(0)\|_{W_\varrho^{2-2/\varrho}(\Omega)}+
\|\theta(0)\|_{W_\varrho^{2-2/\varrho}(\Omega)}\le k_5<\infty$,\\
where $c_0$ is a constant from Lemma 2.2, and $\psi_0$ is an increasing 
function from Lemma 3.3 and $k_1,\dots,k_5$ are given constants.

\proclaim Main Theorem.  
For every fixed $T$, conditions 1--6 with constants $k_1-k_5$, $c_0,c_1$ there 
exist a sufficiently small $d$ and a constant 
$B=B(k_1,\dots,k_5,c_0,c_1)<\infty$ such that for any strong solution 
$v,p,\theta$ to problem (1.1) we have
$$
\|v\|_{W_\varrho^{2,1}(\Omega^t)}+\|\nabla p\|_{L_\varrho(\Omega^t)}+
\|\theta\|_{W_\varrho^{2,1}(\Omega^t)}\le B,
\leqno(1.2)
$$
$$
\|h\|_{W_\sigma^{2,1}(\Omega^t)}+\|\nabla q\|_{L_\sigma(\Omega^t)}+
\|\vartheta\|_{W_\sigma^{2,1}(\Omega^t)}\le B.
\leqno(1.3)
$$
$t\le T$.

\noindent
The result is necessary to prove a long time existence of regular solutions 
to (1.1) in [6].

\section{2. Preliminaries}

In this section we introduce notation and basic estimates for weak solutions 
to problem (1.1).

\subsection{2.1. Notation}

We use isotropic and anisotropic Lebesgue spaces: $L_p(Q)$, 
$Q\in\{\Omega^T,S^T,\Omega,S\}$, $p\in[1,\infty]$; $L_q(0,T;L_p(Q))$, 
$Q\in\{\Omega,S\}$, $p,q\in[1,\infty]$;\\
Sobolev spaces
$$
W_q^{s,s/2}(Q^T),\quad Q\in\{\Omega,S\},\quad q\in[1,\infty],\quad 
s\in\N\cup\{0\}
$$
with the norm
$$
\|u\|_{W_q^{s,s/2}(Q^T)}=\bigg(\sum_{|\alpha|+2a\le s}\intop_{Q^T}
|D_x^\alpha\partial_t^au|^qdxdt\bigg)^{1/q},
$$
where $D_x^\alpha=\partial_{x_1}^{\alpha_1}\partial_{x_2}^{\alpha_2}
\partial_{x_3}^{\alpha_3}$, $|\alpha|=\alpha_1+\alpha_2+\alpha_3$, 
$a,\alpha_i\in\N\cup\{0\}$.

\noindent
In the case $q=2$
$$
H^s(Q)=W_2^s(Q),\quad H^{s,s/2}(Q^T)=W_2^{s,s/2}(Q^T),\quad 
Q\in\{\Omega,S\}.
$$
Moreover, $L_2(Q)=H^0(Q)$, $L_p(Q)=W_p^0(Q)$, $L_p(Q^T)=W_p^{0,0}(Q^T)$.

\noindent
We define a space natural for study weak solutions to the Navier-Stokes and 
parabolic equations
$$
V_2^k(\Omega^T)=\bigg\{u:\ \|u\|_{V_2^k(\Omega^T)}={\rm esssup}_{t\in[0,T]}
\|u\|_{H^k(\Omega)}+\bigg(\intop_0^T\|\nabla u\|_{H^k(\Omega)}^2dt\bigg)^{1/2}
<\infty\bigg\}.
$$

\subsection{2.2. Weak solutions}

By a weak solution to problem (1.1) we mean $v\in V_2^0(\Omega^T)$, 
$\theta\in V_2^0(\Omega^T)\cap L_\infty(\Omega^T)$ satisfying the integral 
identities
$$\eqal{
&-\intop_{\Omega^T}v\cdot\varphi_{,t}dxdt+\intop_{\Omega^T}v\cdot\nabla v\cdot
\varphi dxdt+{\nu\over2}\intop_{\Omega^T}\D(v)\cdot\D(\varphi)dxdt\cr
&=\intop_{\Omega^T}\alpha(\theta)f\cdot\varphi dxdt+
\intop_\Omega v(0)\varphi(0)dx,\cr}
\leqno(2.1)
$$
$$\eqal{
&-\intop_{\Omega^T}\theta\psi_{,t}dxdt+\intop_{\Omega^T}v\cdot\nabla\theta
\psi dxdt+\varkappa\intop_{\Omega^T}\nabla\theta\cdot\nabla\psi dxdt\cr
&=\intop_\Omega\theta(0)\psi(0)dx,\cr}
\leqno(2.2)
$$
which hold for $\varphi,\psi\in W_2^{1,1}(\Omega^T)\cap L_5(\Omega^T)$ 
such that $\varphi(T)=0$, $\psi(T)=0$, $\divv\varphi=0$, 
$\varphi\cdot\bar n|_S=0$.

\proclaim Lemma 2.1. (the Korn inequality, see [10]) 
Assume that
$$
E_\Omega(v)=\|\D(v)\|_{L_2(\Omega)}^2<\infty,\quad 
v\cdot\bar n|_S=0,\quad \divv v=0.
\leqno(2.3)
$$
If $\Omega$ is not axially symmetric there exists a constant $c_1$ such that
$$
\|v\|_{H^1(\Omega)}^2\le c_1E_\Omega(v).
\leqno(2.4)
$$
If $\Omega$ is axially symmetric, $\eta=(-x_2,x_1,0)$, 
$\alpha=\intop_\Omega v\cdot\eta dx$, then there exists a constant $c_2$ 
such that
$$
\|v\|_{H^1(\Omega)}^2\le c_2(E_\Omega(v)+|\alpha|^2).
\leqno(2.5)
$$

\noindent
Let us consider the problem
$$\eqal{
&h_{,t}-\divv\T(h,q)=f\quad &{\rm in}\ \ \Omega^T,\cr
&\divv h=0\quad &{\rm in}\ \ \Omega^T,\cr
&\bar n\cdot h=0,\ \ \bar n\cdot\D(h)\cdot\bar\tau_\alpha=0,\ \ 
\alpha=1,2,\quad &{\rm on}\ \ S_1^T,\cr
&h_i=0,\ \ i=1,2,\ \ h_{3,x_3}=0\quad &{\rm on}\ \ S_2^T,\cr
&h|_{t=0}=h(0)\quad &{\rm in}\ \ \Omega\cr}
\leqno(2.5')
$$

\proclaim Lemma 2.2. 
Let $f\in L_p(\Omega^T)$, $h(0)\in W_p^{2-2/p}(\Omega)$, $S\in C^2$, 
$q<p<\infty$. Then there exists a solution to problem $(2.5')$ such that 
$h\in W_p^{2,1}(\Omega^T)$, $\nabla q\in L_p(\Omega^T)$ and there exists 
a constant $c$ depending on $S$ and $p$ such that
$$\eqal{
&\|h\|_{W_p^{2,1}(\Omega^T)}+\|\nabla q\|_{L_p(\Omega^T)}\le c
(\|f\|_{L_p(\Omega^T)}\cr
&\quad+\|h(0)\|_{W_p^{2-2/p}(\Omega)}).\cr}
\leqno(2.6')
$$

\noindent
The proof is similar to the proof from [1].

\proclaim Lemma 2.3. 
Assume that $v(0)\in L_2(\Omega)$, $\theta(0)\in L_\infty(\Omega)$, 
$f\in L_2(0,T;L_{6/5}(\Omega))$, $T<\infty$. Assume that $\Omega$ is not 
axially symmetric. Assume that there exist constants $\theta_*,\theta^*$ such 
that $\theta_*<\theta^*$ and $\theta_*\le\theta_0(x)\le\theta^*$, 
$x\in\Omega$.\\
Then there exists a weak solution to problem (1.1) such that 
$(v,\theta)\in V_2^0(\Omega^T)\times V_2^0(\Omega^T)$, 
$\theta\in L_\infty(\Omega^T)$ and
$$
\theta_*\le\theta(x,t)\le\theta^*,\quad (x,t)\in\Omega^T,
\leqno(2.6)
$$
$$
\|v\|_{V_2^0(\Omega^T)}\le c(a(\|\theta_0\|_{L_\infty(\Omega)})
\|f\|_{L_2(0,T;L_{6/5}(\Omega))}+\|v_0\|_{L_2(\Omega)})\le c_0,
\leqno(2.7)
$$
$$
\|\theta\|_{V_2^0(\Omega^T)}\le c\|\theta_0\|_{L_2(\Omega)}\le c_0.
\leqno(2.8)
$$

\Proof 
Estimate (2.6) follows from standard considerations (see [5, Lemmas 3.1, 
3.2]). Estimates (2.7), (2.8) follow formally from $(1.1)_{1,3}$ by 
multiplying them by $v$ and $\theta$, respectively, integrating over 
$\Omega$ and $(0,t)$, $t\in(0,T)$, employing (2.6), $(1.1)_2$ and using the 
boundary and initial conditions $(1.1)_{4,5,6,7}$.

\noindent
Existence can be shown in the same way as in [3, Ch. 3, Sect. 1--5].

\noindent
This concludes the proof.
\kwadrat

\Remark{2.4.} 
If $\theta(0)\ge0$, then $\theta(t)\ge0$ for $t\ge0$.

\subsection{2.3. Auxiliary problems}

To prove the existence of global regular solutions we distinguish the 
quantities
$$
h=v_{,x_3},\quad q=p_{,x_3},\quad g=f_{,x_3},\quad \vartheta=\theta_{,x_3}.
\leqno(2.9)
$$
Differentiating $(1.1)_{1,2,4,5}$ with respect to $x_3$ and using [1, 2] yields
$$\eqal{
&h_{,t}-\divv\T(h,q)=-v\cdot\nabla h-h\cdot\nabla v+\alpha_\theta\vartheta f+
\alpha g\quad &{\rm in}\ \ \Omega^T,\cr
&\divv h=0\quad &{\rm in}\ \ \Omega^T,\cr
&\bar n\cdot h=0,\ \ \bar n\cdot\D(h)\cdot\bar\tau_\alpha=0,\ \ \alpha=1,2
\quad &{\rm on}\ \ S_1^T,\cr
&h_i=0,\ \ i=1,2,\ \ h_{3,x_3}=0\quad &{\rm on}\ \ S_2^T,\cr
&h|_{t=0}=h(0)\quad &{\rm in}\ \ \Omega.\cr}
\leqno(2.10)
$$
Let $q$ and $f_3$ be given, Then $w=v_3$ is a solution to the problem
$$\eqal{
&w_{,t}+v\cdot\nabla w-\nu\Delta w=-q+\alpha(\theta)f_3\quad &{\rm in}\ \ 
\Omega^T,\cr
&\bar n\cdot\nabla w=0\quad &{\rm on}\ \ S_1^T,\cr
&w=0\quad &{\rm on}\ \ S_2^T,\cr
&w|_{t=0}=w(0)\quad &{\rm in}\ \ \Omega.\cr}
\leqno(2.11)
$$
Let $F=(\rot f)_3$, $h,v,w$ be given. Then $\chi=(\rot v)_3$ is a solution 
to the problem (see [8])
$$\eqal{
&\chi_{,t}+v\cdot\nabla\chi-h_3\chi+h_2w_{,x_1}-h_1w_{,x_2}-\nu\Delta\chi\quad
\cr
&\quad=\alpha_\theta(\theta_{,x_1}f_2-\theta_{,x_2}f_1)+\alpha F\quad 
&{\rm in}\ \ \Omega^T,\cr
&\chi=v_i(n_{i,x_j}\tau_{1j}+\tau_{1i,x_j}n_j)+v\cdot\bar\tau_1
(\tau_{12,x_1}-\tau_{11,x_2})\cr
&\equiv\chi_*\quad &{\rm on}\ \ S_1^T,\cr
&\chi_{,x_3}=0\quad &{\rm on}\ \ S_2^T,\cr
&\chi|_{t=0}=\chi(0)\quad &{\rm in}\ \ \Omega,\cr}
\leqno(2.12)
$$
where the summation convention over repeated indices is assumed.

\noindent
Differentiating $(1.1)_{3,6,7}$ with respect to $x_3$ yields
$$\eqal{
&\vartheta_{,t}+v\cdot\nabla\vartheta+h\cdot\nabla\theta-
\varkappa\Delta\vartheta=0\quad &{\rm in}\ \ \Omega^T,\cr
&\bar n\cdot\nabla\vartheta=0\quad &{\rm on}\ \ S_1^T,\cr
&\vartheta=0\quad &{\rm on}\ \ S_2^T,\cr
&\vartheta|_{t=0}=\vartheta(0)\quad &{\rm in}\ \ \Omega.\cr}
\leqno(2.13)
$$

\proclaim Lemma 2.5. 
Assume that $\D(h)\in L_2(\Omega)$, $h\cdot\bar n|_S=0$, $\divv h=0$. 
Then $h$ satisfies the inequality
$$
\|h\|_{H^1(\Omega)}\le c\|\D(h)\|_{L_2(\Omega)}.
\leqno(2.14)
$$
where $c$ is a constant independent of $h$.

\Proof 
To show (2.14) we examine the expression
$$
\intop_\Omega|\D(h)|^2dx=\intop_\Omega(h_{i,x_j}+h_{j,x_i})^2dx=
\intop_\Omega(2h_{i,x_j}^2+2h_{i,x_j}h_{j,x_i})dx,
$$
where the second expression implies
$$\eqal{
&\intop_\Omega h_{i,x_j}h_{j,x_i}dx=\intop_\Omega(h_{i,x_j}h_j)_{,x_i}dx
-\intop_\Omega h_{i,x_ix_j}h_jdx=\intop_{S_1\cup S_2}n_ih_{i,x_j}h_jdS\cr
&=-\intop_{S_1}n_{i,x_j}h_ih_jdS_1+\intop_{S_2}n_ih_{i,x_j}h_jdS_2
=-\intop_{S_1}n_{i,x_j}h_ih_jdS_1.\cr}
$$
From the above considerations we have
$$
\|\nabla h\|_{L_2(\Omega)}^2\le c\intop_\Omega|\D(h)|^2dx+c\|h\|_{L_2(S_1)}^2.
\leqno(2.15)
$$
By the trace theorem
$$
\|\nabla h\|_{L_2(\Omega)}^2\le c(\|\D(h)\|_{L_2(\Omega)}^2+
\|h\|_{L_2(\Omega)}^2).
\leqno(2.16)
$$
From [9] we have that
$$
\|h\|_{L_2(\Omega)}\le\delta\|\nabla h\|_{L_2(\Omega)}+
M\|\D(h)\|_{L_2(\Omega)},
\leqno(2.17)
$$
where $\delta$ can be chosen sufficiently small and $M=M(\delta)$ is some 
constant. From (2.15)-\break -(2.17) we have
$$
\|\nabla h\|_{L_2(\Omega)}^2\le c\|\D(h)\|_{L_2(\Omega)}^2.
\leqno(2.18)
$$
From (2.18) and (2.17) we obtain (2.14). This concludes the proof.
\kwadrat

Let us consider the elliptic problem
$$\eqal{
&v_{2,x_1}-v_{1,x_2}=\chi\quad &{\rm in}\ \ \Omega\subset\R^2,
\cr
&v_{1,x_1}+v_{2,x_2}=-h_3\quad &{\rm in}\ \ \Omega\subset\R^2,\cr
&v\cdot\bar n=0\quad &{\rm on}\ \ \partial\Omega.\cr}
\leqno(2.22)
$$
where $x_3$ is treated as a parameter.

\proclaim Lemma 2.6. 
Assume that $\chi,h_3\in L_2(\Omega)$. Then $v\in H^1(\Omega)$ and
$$
\|v\|_{H^1(\Omega)}\le c(\|\chi\|_{L_2(\Omega)}+\|h_3\|_{L_2(\Omega)}).
\leqno(2.23)
$$
Assume that $\chi,h_3\in H^1(\Omega)$. Then $v\in H^2(\Omega)$ and
$$
\|v\|_{H^2(\Omega)}\le c(\|\chi\|_{H^1(\Omega)}+\|h_3\|_{H^1(\Omega)}).
\leqno(2.24)
$$

\Proof
To solve problem (2.22) we introduce the potential $\varphi,\psi$ such that
$$\eqal{
&v_1=\varphi_{,x_1}+\psi_{,x_2},\cr
&v_2=\varphi_{,x_2}-\psi_{,x_1}.\cr}
\leqno(2.25)
$$
Using representation (2.25) we see that $(2.22)_3$ takes the form
$$
\bar n\cdot\nabla\varphi+\bar\tau\cdot\nabla\psi=0\quad {\rm on}\ \ S,
\leqno(2.26)
$$
where $\bar n\perp TS$, $\bar\tau\in TS$.

\noindent
The potential $\varphi$ and $\psi$ are determined up to an arbitrary constant. 
Moreover, to determine the potential we split boundary condition (3.26) into 
two boundary conditions
$$\eqal{
&\bar n\cdot\nabla\varphi|_S=0\cr
&\bar\tau\cdot\nabla\psi|_S=0\Rightarrow\psi|_S=0.\cr}
\leqno(2.27)
$$
Having $v=(v_1,v_2)$ given we calculate $\varphi$ and $\psi$ from the problems
$$\eqal{
&\Delta\varphi=v_{1,x_1}+v_{2,x_2}\quad {\rm in}\ \ \Omega,\cr
&\bar n\cdot\nabla\varphi|_S=0,\cr
&\intop_\Omega\varphi dx=0\cr}
\leqno(2.28)
$$
and
$$\eqal{
&\Delta\psi=v_{1,x_2}-v_{2,x_1}\cr
&\psi|_S=0.\cr}
\leqno(2.29)
$$
In view of (2.28), (2.29) problem (2.22) takes the form
$$\eqal{
&\Delta\psi=\chi\quad &\psi|_S=0,\cr
&\Delta\varphi=-h_3\quad &\bar n\cdot\nabla\varphi|_S=0,\quad 
\intop_\Omega\varphi dx=0\cr}
\leqno(2.30)
$$
Solving problems (2.30) we have the estimates
$$\eqal{
&\|\psi\|_{H^2(\Omega)}\le c\|\chi\|_{L_2(\Omega)},\cr
&\|\varphi\|_{H^2(\Omega)}\le c\|h_3\|_{L_2(\Omega)}.\cr}
\leqno(2.31)
$$
Hence in view of (2.25) we get (2.23).

\noindent
For more regular $\chi$ and $h_3$ we have also the estimates
$$\eqal{
&\|\psi\|_{H^3(\Omega)}\le c\|\chi\|_{H^1(\Omega)},\cr
&\|\varphi\|_{H^3(\Omega)}\le c\|h_3\|_{H^1(\Omega)}.\cr}
\leqno(2.32)
$$
Then (2.25) implies (2.24). This concludes the proof.
\kwadrat

Now we formulate the result on local existence of solutions to problem (1.1) 
with regularity allowed by the regularity of data formulated in the Main 
Theorem. The aim of the result is such that derived in this paper estimates 
are not only a priori type estimates.

\proclaim Lemma 2.7. 
Let the assumptions of the Main Theorem hold.\\
Then for any $A>0$ there exists $t_*>0$ and $(v,\theta,p)$-solution to problem 
(1.1) such that $v\in W_\varrho^{2,1}(\Omega^{t_*})$, 
$\theta\in W_\varrho^{2,1}(\Omega^{t_*})$, 
$\nabla p\in L_\varrho(\Omega^{t_*})$, 
$h\in W_\sigma^{2,1}(\Omega^{t_*})$, $\nabla q\in L_\sigma(\Omega^{t_*})$ and
$$\eqal{
&\|h\|_{W_\sigma^{2,1}(\Omega^{t_*})}+\|\nabla q\|_{L_\sigma(\Omega^{t_*})}
\le A,\cr
&\|v\|_{W_\varrho^{2,1}(\Omega^{t_*})}+
\|\theta\|_{W_\varrho^{2,1}(\Omega^{t_*})}+
\|\nabla p\|_{L_\varrho(\Omega^{t_*})}\le A.\cr}
$$

\noindent
Consider the problem
$$\eqal{
&u_{,t}-\nu\Delta u=0\cr
&u|_S=\varphi\cr
&u|_{t=0}=0\cr}
$$

\proclaim Lemma 2.8. (see [4]) 
Let $\varphi\in L_q(0,T;L_p(S))$, $p,q\in[1,\infty]$ then 
$u\in L_q(0,T;L_p(\Omega))$ and
$$
\|u\|_{L_q(0,T;L_p(\Omega))}\le c\|\varphi\|_{L_q(0,T;L_p(S))}.
$$
Let $\varphi\in W_2^{{1\over2},{1\over4}}(S^T)$ then 
$u\in W_2^{1,1/2}(\Omega^T)$ and
$$
\|u\|_{W_2^{1,1/2}(\Omega^T)}\le c\|\varphi\|_{W_2^{1/2,1/4}(S^T)}.
$$

\section{3. Estimates}

\proclaim Lemma 3.1. 
Let the assumptions of Lemma 2.3 be satisfied. Moreover assume that 
$f\in L_2(0,T;L_3(\Omega))$, $f_3\in L_3(0,T;L_{4\over3}(S_2))$, 
$g\in L_2(0,T;L_{6/5}(\Omega))$, $h(0)\in L_2(\Omega)$, 
$\vartheta(0)\in L_2(\Omega)$, $\nabla v\in L_2(0,T;L_3(\Omega))$, 
$\nabla\theta\in L_2(0,T;L_3(\Omega))$. Assume that solutions to (2.10), 
(2.13) are sufficiently regular. Let $c_1=a(\|\theta_0\|_{L_\infty})$.\\
Let $h\in(0,T;L_3(\Omega))$, then solutions of (2.10), (2.13) satisfy the 
inequality
$$\eqal{
&\|h\|_{V_2^0(\Omega^T)}+\|\vartheta\|_{V_2^0(\Omega^t)}^2\le c\exp
(cc_1^2\|f\|_{L_2(0,t;L_3(\Omega))}^2)\cr
&\quad\cdot[c_0^2\|h\|_{L_\infty(0,t;L_3(\Omega))}^2+c_1^2
\|g\|_{L_2(0,t;L_{6/5}(\Omega))}^2+c_1^2\|f_3\|_{L_2(0,t;L_{4/3}(S_2))}^2\cr
&\quad+\|h(0)\|_{L_2(\Omega)}^2
+\|\vartheta(0)\\_{L_2(\Omega)}^2],\quad t\le T.\cr}
\leqno(3.1)
$$
Let, additionally, $v,\theta\in L_2(0,T;W_3^1(\Omega))$, then for solutions 
of (2.10), (2.13) we have
$$\eqal{
&\|h\|_{V_2^0(\Omega^t)}^2+\|\vartheta\|_{V_2^0(\Omega^t)}^2\le c\exp
[c(\|\nabla v\|_{L_2(0,t;L_3(\Omega))}^2\cr
&\quad+\|\nabla\theta\|_{L_2(0,t;L_3(\Omega))}^2+c_1^2
\|f\|_{L_2(0,t;L_3(\Omega))}^2)]\cr
&\quad\cdot[c_1^2\|g\|_{L_2(0,t;L_{6/5}(\Omega))}^2+c_1^2
\|f_3\|_{L_2(0,t;L_{4\over3}(S_2))}^2+\|h(0)\|_{L_2(\Omega)}^2+
\|\vartheta(0)\|_{L_2(\Omega)}^2]\cr}
\leqno(3.2)
$$
$t\le T$.

\Proof 
Multiplying (2.10) by $h$, integrating over $\Omega$ and using Lemma 2.5 yields
$$\eqal{
&{1\over2}{d\over dt}\|h\|_{L_2(\Omega)}^2+\nu\|h\|_{H^1(\Omega)}^2\le c
\intop_\Omega|h\cdot\nabla v\cdot h|dx+c\intop_\Omega|\alpha_\theta\vartheta fh|dx\cr
&\quad+c\intop_\Omega|\alpha gh|dx+c\intop_{S_2}|\alpha f_3h_3|dx_1dx_2\cr}
\leqno(3.3)
$$
where the first term on the r.h.s. we estimate by
$$
\varepsilon_1\|h\|_{L_6(\Omega)}^2+c(1/\varepsilon_1)
\|\nabla v\|_{L_2(\Omega)}^2\|h\|_{L_3(\Omega)}^2,
$$
the second by
$$
\varepsilon_2\|h\|_{L_6(\Omega)}^2+c(1/\varepsilon_2)a^2
(\|\theta_0\|_{L_\infty(\Omega)})\|\vartheta f\|_{L_{6/5}(\Omega)}^2,
$$
the third by
$$
\varepsilon_3\|h\|_{L_6(\Omega)}^2+c(1/\varepsilon_3)a^2
(\|\theta_0\|_{L_\infty(\Omega)})\|g\|_{L_{6/5}(\Omega)}^2
$$
and the fourth by
$$
\varepsilon_4\|h\|_{H^1(\Omega)}^2+c(1/\varepsilon_4)a^2
(\|\theta_0\|_{L_\infty(\Omega)})\|f_3\|_{L_{4\over3}(S_2)}.
$$
Assuming that $\varepsilon_1,\varepsilon_2,\varepsilon_3,\varepsilon_4$ are 
sufficiently small we obtain
$$\eqal{
&{d\over dt}\|h\|_{L_2(\Omega)}^2+\nu\|h\|_{H^1(\Omega)}^2\le c
(\|\nabla v\|_{L_2(\Omega)}^2\|h\|_{L_3(\Omega)}^2\cr
&\quad+c_1^2(\|\vartheta\|_{L_2(\Omega)}^2\|f\|_{L_3(\Omega)}^2+
\|g\|_{L_{6/5}(\Omega)}^2+\|f_3\|_{L_{4\over3}(S_2)}^2)).\cr}
\leqno(3.4)
$$
Multiplying (2.13) by $\vartheta$ and integrating over $\Omega$ yields
$$\eqal{
&{1\over2}{d\over dt}\|\vartheta\|_{L_2(\Omega)}^2+\varkappa
\|\vartheta\|_{H^1(\Omega)}^2\le c\intop_\Omega|h\cdot\nabla\theta\vartheta|dx
\cr
&\le\varepsilon\|\vartheta\|_{L_6(\Omega)}^2+c(1/\varepsilon)
\|h\|_{L_3(\Omega)}^2\|\nabla\theta\|_{L_2(\Omega)}^2.\cr}
\leqno(3.5)
$$
For sufficiently small $\varepsilon$ we have
$$
{d\over dt}\|\vartheta\|_{L_2(\Omega)}^2+\varkappa
\|\vartheta\|_{H^1(\Omega)}^2\le c\|h\|_{L_3(\Omega)}^2
\|\nabla\theta\|_{L_2(\Omega)}^2.
\leqno(3.6)
$$
Adding (3.4) and (3.6), integrating with respect to time and using (2.7) and 
(2.8) we obtain (3.1).

We can replace inequalities (3.4) and (3.6) by
$$\eqal{
&{d\over dt}\|h\|_{L_2(\Omega)}^2+\nu\|h\|_{H^1(\Omega)}^2\le c
(\|\nabla v\|_{L_3(\Omega)}^2\|h\|_{L_2(\Omega)}^2\cr
&\quad+c_1^2(\|\vartheta\|_{L_2(\Omega)}^2\|f\|_{L_3(\Omega)}^2+
\|g\|_{L_{6/5}(\Omega)}^2+\|f_3\|_{L_{4\over3}(S_2)}^2))\cr}
\leqno(3.7)
$$
and
$$
{d\over dt}\|\vartheta\|_{L_2(\Omega)}^2+\varkappa
\|\vartheta\|_{H^1(\Omega)}^2\le c\|\nabla\theta\|_{L_3(\Omega)}^2
\|h\|_{L_2(\Omega)}^2.
\leqno(3.8)
$$
Adding (3.7) and (3.8), integrating the sum with respect to time yields (3.2). 
This ends the proof.
\kwadrat

\noindent
To obtain an estimate for solutions to problem (2.12) we introduce a function 
$\tilde\chi:\Omega\times[0,T]\to\R$ as a solution to the problem
$$\eqal{
&\tilde\chi_{,t}-\nu\Delta\tilde\chi=0\quad &{\rm in}\ \ \Omega^T,\cr
&\tilde\chi=\chi_*\quad &{\rm on}\ \ S_1^T,\cr
&\tilde\chi_{,x_3}=0\quad &{\rm on}\ \ S_2^T,\cr
&\tilde\chi|_{t=0}=0\quad &{\rm in}\ \ \Omega.\cr}
\leqno(3.9)
$$
Then the function
$$
\chi'=\chi-\tilde\chi
\leqno(3.10)
$$
satisfies
$$\eqal{
&\chi'_{,t}+v\cdot\nabla\chi'-h_3\chi'+h_2w_{,x_1}-h_1w_{,x_2}-\nu\Delta\chi'
\quad\cr
&=\alpha_\theta(\theta_{,x_1}f_2-\theta_{,x_2}f_1)+\alpha F-
v\cdot\nabla\tilde\chi+h_3\tilde\chi\quad &{\rm in}\ \ \Omega^T,\cr
&\chi'=0\quad &{\rm on}\ \ S_1^T,\cr
&\chi'_{,x_3}=0\quad &{\rm on}\ \ S_2^T,\cr
&\chi'|_{t=0}=\chi(0)\quad &{\rm in}\ \ \Omega.\cr}
\leqno(3.11)
$$

\proclaim Lemma 3.2. 
Let the assumptions of Lemma 2.3 be satisfied. Moreover assume that 
$h,f\in L_\infty(0,T;L_3(\Omega))$, $F\in L_2(0,T;L_{6/5}(\Omega))$, 
$v'=(v_1,v_2)\in L_\infty(0,T;H^{1/2+\varepsilon}(\Omega))\cap 
W_2^{1,1/2}(\Omega^T)$, $\chi(0)\in L_2(\Omega)$, $\varepsilon>0$ 
is arbitrary small.\\
Assume that solutions to (1.1) are sufficiently regular. 
Then for solutions to (2.12) we have
$$\eqal{
&\|\chi\|_{V_2^0(\Omega^t)}^2\le c\bigg(c_0^2\sup_t\|h\|_{L_3(\Omega)}^2+
c_1^2c_0^2\sup_t\|f\|_{L_3(\Omega)}^2\cr
&\quad+c_1^2\|F\|_{L_2(0,t;L_{6/5}(\Omega))}^2+c_0^2\varepsilon_7^2
\|v'\|_{L_\infty(0,t;H^1(\Omega))}^2\cr
&\quad+\|v'\|_{L_\infty(0,t;H^{1/2+\varepsilon}(\Omega))}^2+
\|v'\|_{W_2^{1,1/2}(\Omega^t)}+\|\chi(0)\|_{L_2(\Omega)}^2\cr
&\quad+\bigg(c_0^2c^2\bigg({1\over\varepsilon_7}\bigg)+
\sup_t\|h\|_{L_3(\Omega)}^2\bigg)(a^2(\|\theta_0\|_{L_\infty(\Omega^t)})
\|f\|_{L_2(0,t;L_{6/5}(\Omega))}^2+\|v_0\|_{L_2(\Omega)}^2)\bigg),\cr
&t\le T.\cr}
\leqno(3.12)
$$

\Proof 
Multiplying $(3.11)_1$ by $\chi'$, integrating over $\Omega$, using boundary 
conditions $(3.11)_{2,3}$, $(1.1)_5$ and $(1.1)_2$, we obtain
$$\eqal{
&{1\over2}{d\over dt}\|\chi'\|_{L_2(\Omega)}^2+\nu
\|\nabla\chi'\|_{L_2(\Omega)}^2=\intop_\Omega h_3\chi'^2dx\cr
&\quad-\intop_\Omega(h_2w_{,x_1}-h_1w_{,x_2})\chi'dx+
\intop_\Omega\alpha_\theta(\theta_{,x_1}f_2-\theta_{,x_2}f_1)\chi'dx\cr
&\quad+\intop_\Omega\alpha F\chi'dx-\intop_\Omega v\cdot\nabla\tilde\chi
\chi'dx+\intop_\Omega h_3\tilde\chi\chi'dx.\cr}
\leqno(3.13)
$$
Now we estimate the terms on the r.h.s. of the above equality. Let 
$x'=(x_1,x_2)$. The first term we estimate by
$$
\bigg|\intop_\Omega h_3\chi'^2dx\bigg|\le
\varepsilon_1\|\chi'\|_{L_6(\Omega)}^2+{c\over\varepsilon_1}
\|\chi'\|_{L_2(\Omega)}^2\|h_3\|_{L_3(\Omega)}^2,
$$
the second by
$$
\varepsilon_2\|\chi'\|_{L_6(\Omega)}^2+{c\over\varepsilon_2}
\|h\|_{L_3(\Omega)}^2\|w_{,x'}\|_{L_2(\Omega)}^2,
$$
the third by
$$
\varepsilon_3\|\chi'\|_{L_6(\Omega)}^2+{c\over\varepsilon_3}c_1^2
\|\theta_{,x}\|_{L_2(\Omega)}^2\|f\|_{L_3(\Omega)}^2,
$$
where we used $(1.1')$. The fourth by
$$
\varepsilon_4\|\chi'\|_{L_6(\Omega)}^2+{c\over\varepsilon_4}c_1^2
\|F\|_{L_{6/5}(\Omega)}^2,
$$
where we used also $(1.1')$.

\noindent
To estimate the fifth term on the r.h.s. of (3.13) we integrate it by parts 
and use $(1.1)_{2,5}$. Then it takes the form
$$
I\equiv\intop_\Omega v\cdot\nabla\chi'\tilde\chi dx.
$$
Hence
$$
|I|\le\varepsilon_5\|\nabla\chi'\|_{L_2(\Omega)}^2+{c\over\varepsilon_5}
\|v\|_{L_6(\Omega)}^2\|\tilde\chi\|_{L_3(\Omega)}^2.
$$
Finally, the last term on the r.h.s. of (3.13) is bounded by
$$
\varepsilon_6\|\chi'\|_{L_6(\Omega)}^2+{c\over\varepsilon_6}
\|h\|_{L_3(\Omega)}^2\|\tilde\chi\\_{L_2(\Omega)}^2.
$$
Using the above estimates in (3.13), assuming that 
$\varepsilon_1,\dots,\varepsilon_6$ are sufficiently small, integrating 
the result with respect to time and using (2.7)--(2.8) we obtain
$$\eqal{
&\|\chi'\|_{V_2^0(\Omega^t)}^2\le c(\sup_t\|h\|_{L_3(\Omega)}^2
\|\chi'\|_{L_2(0,t;L_2(\Omega))}^2\cr
&\quad+c_0^2\sup_t\|h\|_{L_3(\Omega)}^2+c_1^2c_0^2\sup_t
\|f\|_{L_3(\Omega)}^2+c_1^2\|F\|_{L_2(0,t;L_{6/5}(\Omega))}^2\cr
&\quad+c_0^2\|\tilde\chi\|_{L_\infty(0,t;L_3(\Omega))}^2+
\sup_t\|h\|_{L_3(\Omega)}^2\|\tilde\chi\|_{L_2(0,t;L_2(\Omega))}^2+
\|\chi(0)\|_{L_2(\Omega)}^2).\cr}
\leqno(3.14)
$$
In view of (2.7) we have $\|\chi\|_{L_2(\Omega^t)}\le cc_0$.

\noindent
Using (3.10) and this fact we obtain from (3.14) the inequality
$$\eqal{
&\|\chi\|_{V_2^0(\Omega^t)}^2\le c(c_0^2\sup_t\|h\|_{L_3(\Omega)}^2+
c_1^2c_0^2\sup_t\|f\|_{L_3(\Omega)}^2\cr
&\quad+c_1^2\|F\|_{L_2(0,t;L_{6/5}(\Omega))}^2+c_0^2
\|\tilde\chi\|_{L_\infty(0,t;L_3(\Omega))}^2+\sup_t\|h\|_{L_3(\Omega)}^2
\|\tilde\chi\|_{L_2(\Omega^t)}^2\cr
&\quad+\|\tilde\chi\|_{V_2^0(\Omega^t)}+\|\chi(0)\|_{L_2(\Omega)}^2).\cr}
\leqno(3.15)
$$
Since $\tilde\chi$ is a solution of (3.9) and $\chi_*$ is described by 
$(2.12)_2$ we have the estimates by Lemma 2.8,
$$\eqal{
&\intop_0^t\|\tilde\chi(t')\|_{L_2(\Omega)}^2dt'\le c\intop_0^t
\|v'(t')\|_{L_2(S)}^2dt'\le c\intop_0^t\|v'(t')\|_{H^1(\Omega)}^2dt'\cr
&\le c(a^2(\|\theta_0\|_{L_\infty(\Omega)})
\|f\|_{L_2(0,t;L_{6/5}(\Omega))}^2+\|v_0\|_{L_2(\Omega)}^2),\cr
&\|\tilde\chi\|_{L_\infty(0,t;L_3(\Omega))}\le c\|v'\|_{L_\infty(0,t;L_3(S))}
\le\varepsilon_7\|v'\|_{L_\infty(0,t;H^1(\Omega))}+
c\bigg({1\over\varepsilon_7}\bigg)\|v'\|_{L_\infty(0,t;L_2(\Omega))},\cr
&\quad\|\tilde\chi\|_{V_2^0(\Omega^t)}^2\le c
\bigg(\|\tilde\chi\|_{L_\infty(0,t;L_2(\Omega))}^2+\intop_0^t
\|\tilde\chi(t')\|_{H^1(\Omega)}^2dt'\bigg)\cr
&\le c\|v'\|_{L_\infty(0,t;H^{1/2+\varepsilon}(\Omega))}^2+c
\|v'\|_{W_2^{1,1/2}(\Omega^t)}^2,\quad \varepsilon>0.\cr}
\leqno(3.16)
$$
Employing (3.16) in (3.15) yields (3.12). This concludes the proof.
\kwadrat

\noindent
Let us consider the problem
$$\eqal{
&v_{2,x_1}-v_{1,x_2}=\chi\quad &{\rm in}\ \ \Omega',\cr
&v_{1,x_1}+v_{2,x_2}=-h_3\quad &{\rm in}\ \ \Omega',\cr
&v'\cdot\bar n'=0\quad &{\rm on}\ \ S',\cr}
\leqno(3.17)
$$
where $\Omega'=\Omega\cap\{{\rm plane}\ x_3=\const\in(-a,a)\}$, 
$S'=S\cap\{{\rm plane}\ x_3=\const\in(-a,a)\}$, $x_3,t$ are treated as 
parameters, $\bar n'=(n_1,n_2)$.

\proclaim Lemma 3.3. 
Let the assumptions of Lemmas 2.3, 3.1, 3.2 be satisfied. Assume that 
$(v,p,\theta)$ is a weak solution to problem (1.1). Assume that
$$\eqal{
&c_1\|g\|_{L_2(0,t;L_{6/5}(\Omega))}+c_1c_0\|f\|_{L_\infty(0,t;L_3(\Omega))}+
c_1\|F\|_{L_2(0,t;L_{6/5}(\Omega))}\cr
&\quad+c_1\|f_3\|_{L_2(0,t;L_{4\over3}(S_2))}+\|h(0)\|_{L_2(\Omega)}+
\|\vartheta(0)\|_{L_2(\Omega)}+\|\chi(0)\|_{L_2(\Omega)}\cr
&\quad+c_0^2(c_1\|f\|_{L_2(0,t;L_{6/5}(\Omega))}+\|v(0)\|_{L_2(\Omega)})
\le k_1<\infty,\cr
&\|f\|_{L_2(0,t;L_3(\Omega))}\le k_2<\infty,\cr}
\leqno(3.18)
$$
$t\le T$. Then the following inequality
$$
\|v'\|_{V_2^1(\Omega^t)}^2\le c[e^{cc_1^2k_2^2}(c_0^2
\|h\|_{L_\infty(0,t,L_3(\Omega))}^2+\psi(c_0)k_1^2)+
\|v'\|_{L_2(\Omega;H^{1/2}(0,t))}^2]
\leqno(3.19)
$$
holds, where $v'=(v_1,v_2)$, $t\le T$ and $\psi$ is an increasing 
positive function.

\Proof 
Assuming that $\varepsilon_7$ is sufficiently small in view of (3.1), (3.12) 
and Lemma 2.6 we obtain for solutions to problem (3.17) the inequality 
(see [9])
$$\eqal{
\|v'\|_{L_{10}(\Omega^T)}^2&\le c
\|v'\|_{V_2^1(\Omega^t)}^2\le c[e^{cc_1^2k_2^2}(c_0^2
\|h\|_{L_\infty(0,t;L_3(\Omega))}^2+\psi(c_0)k_1^2)\cr
&\quad+\|v'\|_{L_\infty(0,t;H^{1/2+\varepsilon}(\Omega))}^2+
\|v'\|_{W_2^{1,1/2}(\Omega^t)}^2],\cr}
\leqno(3.20)
$$
where $\varepsilon$ is arbitrary small number.

\noindent
By interpolation inequalities
$$
\|v'\|_{L_\infty(0,t;H^{1/2+\varepsilon}(\Omega))}\le\varepsilon_1
\|v'\|_{L_\infty(0,t;H^1(\Omega))}+c(1/\varepsilon_1)
\|v'\|_{L_\infty(0,t;L_2(\Omega))},
$$
and
$$
\|v'\|_{W_2^{1,1/2}(\Omega^t)}=\|v'\|_{L_2(0,t;H^1(\Omega))}+
\|v'\|_{L_2(\Omega;H^{1/2}(0,t))},
$$
where
$$
\|v'\|_{L_2(0,t;H^1(\Omega))}\le ck_1.
$$
Then we obtain (3.19) from (3.20) for sufficiently small $\varepsilon_1$. This 
concludes the proof.
\kwadrat

\noindent
Let us consider problem $(1.1)_{1,2,4,5,7}$ in the form
$$\eqal{
&v_{,t}-\divv\T(v,p)=-v'\cdot\nabla'v-wh+\alpha(\theta)f\quad 
&{\rm in}\ \ \Omega^T,\cr
&\divv v=0\quad &{\rm in}\ \ \Omega^T,\cr
&v\cdot\bar n=0,\ \ \bar n\cdot\D(v)\cdot\bar\tau_\alpha=0,\ \ \alpha=1,2,
\quad &{\rm on}\ \ S^T,\cr
&v|_{t=0}=v_0\quad &{\rm in}\ \ \Omega.\cr}
\leqno(3.21)
$$
where $v'\cdot\nabla'=v_1\partial_{x_1}+v_2\partial_{x_2}$.

\proclaim Lemma 3.4. 
Assume that $(v,p,\theta)$ is a weak solution to problem (1.1). Let the 
assumptions of Lemma 3.3 be satisfied. Let
$$
\|f\|_{L_2(\Omega^t)}+\|v_0\|_{H^1(\Omega)}\le k_3<\infty
$$
$$
H(t)=\|h\|_{L_\infty(0,t;L_3(\Omega))}+\|h\|_{L_{10\over3}(\Omega^t)}<\infty,
$$
$t\le T$.
Then there exists a constant $c_2=c_2(c_0,c_1)$ such that for solutions to 
problem (3.21) the inequality
$$
\|v\|_{W_2^{2,1}(\Omega^t)}+\|\nabla p\|_{L_2(\Omega^t)}\le c_2e^{cc_1^2k_2^2}
(H+1+k_1+k_3)^2+ck_3,\quad t\le T,
\leqno(3.22)
$$
holds.

\noindent
Proof is the same as the proof of Lemma 3.3 in [7].

\noindent
Finally, we obtain an estimate for $h$.

\proclaim Lemma 3.5. 
Let the assumptions of Lemma 3.4 be satisfied. Let
$$\eqal{
&c_1\|f\|_{L_\infty(\Omega^t)}e^{cc_1^2k_2^2}k_1+c_1\|g\|_{L_\sigma(\Omega^t)}
+\|\vartheta(0)\|_{W_\sigma^{2-2/\sigma}(\Omega^t)}+
\|h(0)\|_{W_\sigma^{2-1/\sigma}(\Omega)}\le k_4<\infty,\cr
&c_1\|g\|_{L_2(0,t;L_{6/5}(\Omega))}+c_1\|f_3\|_{L_2(0,t;L_{4\over3}(S_2))}+
\|h(0)\|_{L_2(\Omega)}+\|\vartheta(0)\|_{L_2(\Omega)}\le d<\infty,\cr
&c_1\|f\|_{L_\varrho(\Omega^T)}+\|v(0)\|_{W_\varrho^{2-2/\varrho}(\Omega)}+
\|\theta(0)\|_{W_\varrho^{2-2/\varrho}(\Omega)}\le k_5<\infty,\cr}
\leqno(3.23)
$$
for $t\le T$.
Then for $d$ sufficiently small there exists a constant $A$ such that
$$
\|h\|_{W_\sigma^{2,1}(\Omega^t)}+\|\nabla q\|_{L_\sigma(\Omega^t)}\le A,\quad
{5\over3}<\sigma,\quad t\le T,
\leqno(3.24)
$$
$$
\|\nabla p\|_{L_\varrho(\Omega^t)}+\|v\|_{W_\varrho^{2,1}(\Omega^t)}+
\|\theta\|_{W_\varrho^{2,1}(\Omega^t)}\le\varphi(A)+ck_5,\quad
{5\over3}\le\varrho,\quad t\le T,
\leqno(3.25)
$$
where $\varphi$ is some positive increasing function.

\Proof 
In view of Lemma 2.2 for solutions to problem (2.10) we have
$$\eqal{
&\|h\|_{W_\sigma^{2,1}(\Omega^t)}+\|\nabla q\|_{L_\sigma(\Omega^t)}\le c
(\|v\cdot\nabla h\|_{L_\sigma(\Omega^t)}\cr
&\quad+\|h\cdot\nabla v\|_{L_\sigma(\Omega^t)}+
\|\alpha_\theta\vartheta f\|_{L_\sigma(\Omega^t)}+
\|\alpha g\|_{L_\sigma(\Omega^t)}\cr
&\quad+\|h(0)\|_{W_\sigma^{2-2/\sigma}(\Omega)}).\cr}
\leqno(3.26)
$$
In view of the imbedding
$$
\|v\|_{L_{10}(\Omega^t)}+\|\nabla v\|_{L_{10\over3}(\Omega^t)}\le c
\|v\|_{W_2^{2,1}(\Omega^t)}.
\leqno(3.27)
$$
and inequality (3.22) we estimate the first term on the r.h.s. of (3.26) by
$$
\|v\|_{L_{10}(\Omega^t)}(\varepsilon_1\|h\|_{W_\sigma^{2,1}(\Omega^t)}+
c\bigg({1\over\varepsilon_1}\bigg)\|h\|_{L_2(\Omega^t)})
$$
and the second by
$$
\|\nabla v\|_{L_{10\over3}(\Omega^t)}(\varepsilon_2
\|h\|_{W_\sigma^{2,1}(\Omega^t)}+c\bigg({1\over\varepsilon_2}\bigg)
\|h\|_{L_2(\Omega^t)}).
$$
In view of (2.6) and $(1.1')$ the third and the fourth terms on the r.h.s. 
of (3.26) can be estimated by
$$
cc_1(\|f\|_{L_\infty(\Omega^t)}\|\vartheta\|_{L_\sigma(\Omega^t)}+
\|g\|_{L_\sigma(\Omega^t)})\equiv I.
$$
We use (3.1) with notation (3.18). Then we obtain
$$
I\le cc_1(\|f\|_{L_\infty(\Omega^t)}e^{cc_1^2k_2^2}(k_1+c_0
\|h\|_{L_\infty(0,t;L_3(\Omega))})+\|g\|_{L_\sigma(\Omega^t)}).
$$
We will use also the interpolation
$$
\|h\|_{L_\infty(0,t;L_3(\Omega))}\le\varepsilon_2
\|h\|_{W_\sigma^{2,1}(\Omega^t)}+c(1/\varepsilon_3)\|h\|_{L_2(\Omega^t)}.
$$
Employing the above estimates in (3.26), assuming that 
$\varepsilon_1,\varepsilon_2,\varepsilon_3$ are sufficiently small and 
using (3.22) we obtain
$$\eqal{
&\|h\|_{W_\sigma^{2,1}(\Omega^t)}+\|\nabla q\|_{L_\sigma(\Omega^t)}\le
\varphi(H)\|h\|_{L_2(\Omega^t)}\cr
&\quad+cc_1(\|f\|_{L_\infty(\Omega^t)}e^{cc_1^2k_2^2}k_1+
\|g\|_{L_\sigma(\Omega^t)})+c\|h(0)\|_{W_\sigma^{2-2/\sigma}(\Omega)},\cr}
\leqno(3.28)
$$
where $\varphi$ is an increasing positive function depending on $H$ and 
on constants\\ $c_0,c_1,k_1,\dots,k_5$.

\noindent
Using notation $(3.23)_1$ we have
$$
\|h\|_{W_\sigma^{2,1}(\Omega^t)}+\|\nabla q\|_{L_\sigma(\Omega^t)}\le
\varphi(H)\|h\|_{L_2(\Omega^t)}+ck_4.
\leqno(3.29)
$$
We want to estimate $\|h\|_{L_2(\Omega^t)}$ by applying (3.2). 
For this purpose we need to estimate 
$\|\nabla\theta\|_{L_2(0,t;L_3(\Omega))}$. Hence we consider problem 
$(1.1)_{3,6,7}$ and we are looking for solutions of this problem such that 
$\theta\in W_\varrho^{2,1}(\Omega^t)$ with so large $\varrho$ that
$$
\|\nabla\theta\|_{L_2(0,t;L_3(\Omega))}\le c
\|\theta\|_{W_\varrho^{2,1}(\Omega^t)}.
\leqno(3.30)
$$
We see that (3.30) holds for $\varrho\ge{5\over3}$.

\noindent
Considering problem $(1.1)_{3,6,7}$ we have
$$
\|\theta\|_{W_\varrho^{2,1}(\Omega^t)}\le c
\|v\cdot\nabla\theta\|_{L_\varrho(\Omega^t)}+
\|\theta(0)\|_{W_\varrho^{2-2/\varrho}(\Omega)}).
\leqno(3.31)
$$
The first term on the r.h.s. we estimate by
$$
\|v\|_{L_{\varrho\lambda_1}(\Omega^t)}
\|\nabla\theta\|_{L_{\varrho\lambda_2}(\Omega^t)}\equiv I_1,
$$
where $1/\lambda_1+1/\lambda_2=1$, $\varrho\lambda_1=10$.

\noindent
Using the interpolation inequality
$$
\|\nabla\theta\|_{L_{\varrho\lambda_2}(\Omega^t)}\le\varepsilon_4
\|\theta\|_{W_\varrho^{2,1}(\Omega^t)}+c\bigg({1\over\varepsilon_4}\bigg)
\|\theta\|_{L_2(\Omega^t)}
$$
which holds for ${5\over\varrho}-{5\over\varrho\lambda_2}<1$ so for 
${5\over\varrho\lambda_1}<1$. Hence
$$
I_1\le\|v\|_{L_{10}(\Omega^t)}(\varepsilon_4
\|\theta\|_{W_\varrho^{2,1}(\Omega^t)}+c\bigg({1\over\varepsilon_4}\bigg)
\|\theta\|_{L_2(\Omega^t)}).
$$
Using the estimate in (3.31), assuming that $\varepsilon_4$ is sufficiently 
small, using (3.27) and (3.22), we obtain
$$
\|\theta\|_{W_\varrho^{2,1}(\Omega^t)}\le\varphi(H)+c
\|\theta(0)\|_{W_\varrho^{2-2/\varrho}(\Omega^T)}.
\leqno(3.32)
$$
Similarly by [4, Theorem 2.2]
$$
\|v\|_{W_\varrho^{2,1}(\Omega^t)}+\|\nabla p\|_{L_\varrho(\Omega^t)}\le
\varphi(H)+c_1\|f\|_{L_\varrho(\Omega^T)}+c
\|v(0)\|_{W_\varrho^{2-2/\varrho}(\Omega)}
\leqno(3.33)
$$
Let us consider (3.29). In view of (3.2) we estimate the norm 
$\|h\|_{L_2(\Omega^t)}$, where
$$
\|\nabla v\|_{L_2(0,t;L_3(\Omega))}+\|\nabla\theta\|_{L_2(0,t;L_3(\Omega))}
\le\varphi(H)+ck_5
$$
Then (3.29) takes the form
$$
\|h\|_{W_\sigma^{2,1}(\Omega^t)}+\|\nabla q\|_{L_\sigma(\Omega^t)}\le
\varphi(H)d+ck_4,
\leqno(3.34)
$$
where $\varphi$ is an increasing positive function.

\noindent
Let $\sigma$ be such that
$$
H=\|h\|_{L_\infty(0,t;L_3(\Omega))}+\|h\|_{L_{10\over3}(\Omega^t)}\le c
\|h\|_{W_\sigma^{2,1}(\Omega^t)},
$$
which holds for $\sigma>{5\over3}$.

\noindent
Then (3.34) takes the form
$$
\|h\|_{W_\sigma^{2,1}(\Omega^t)}+\|\nabla q\|_{L_\sigma(\Omega^t)}\le
\varphi(\|h\|_{W_\sigma^{2,1}(\Omega^t)})d+ck_4
\leqno(3.35)
$$
Hence for $d$ sufficiently small there exists a constant $A$ such that
$$
\|h\|_{W_\sigma^{2,1}(\Omega^t)}+\|\nabla q\|_{L_\sigma(\Omega^t)}\le A,
\quad t\le T.
\leqno(3.36)
$$
By (3.36), (3.32) and (3.33) the proof is complete.
\kwadrat
\vskip6pt

\noindent
{\bf Proof of the main Theorem}

Now we want to increase regularity described by (3.25). Assume 
$10\le\varrho<\infty$. In view of [5, Theorem 2.1] for a solution $v$ 
to problem (1.1) we have
$$\eqal{
&\|v\|_{W_\varrho^{2,1}(\Omega^t)}+\|\nabla p\|_{L_\varrho(\Omega^t)}\le c
(\|v\cdot\nabla v\|_{L_\varrho(\Omega^t)}\cr
&\quad+\|\alpha(\theta)f\|_{L_\varrho(\Omega^t)}+
\|v_0\|_{W_\varrho^{2-2/\varrho}(\Omega)}).\cr}
\leqno(3.37)
$$
We estimate the first term on the r.h.s. of (3.37) by
$$\eqal{
&\|v\|_{L_\infty(\Omega^t)}\|\nabla v\|_{L_\varrho(\Omega^t)}\cr
&\le c\|v\|_{W_5^{2,1}(\Omega^t)}(\varepsilon_1
\|v\|_{W_\varrho^{2,1}(\Omega^t)}+c(1/\varepsilon_1)\|v\|_{L_2(\Omega^t)})\cr}
\leqno(3.38)
$$
and the second by
$$
c_1\|f\|_{L_\infty(\Omega^t)}.
\leqno(3.39)
$$
Assuming that $\varepsilon_1$ is sufficiently small and using (3.37)--(3.39) 
we obtain
$$
\|v\|_{W_\varrho^{2,1}(\Omega^t)}+\|\nabla p\|_{L_\varrho(\Omega^t)}\le B_1,
\leqno(3.40)
$$
where $B_1$ is a constant depending on constants from imbedding theorems 
and data.

\noindent
Similarly by [3, Ch. 4, Sect. 9, Th. 9.1] we obtain
$$
\|\theta\|_{W_\varrho^{2,1}(\Omega^t)}\le B_2.
\leqno(3.41)
$$
Now we want to increase regularity described by (3.24).

\noindent
There exist $p'>\sigma$, $p''>{5\over2}$ such that
$$
{5\over\varrho}-{5\over p'}<1,\quad 
{5\over\varrho}-{5\over p''}<1.
$$
Hence $p=\max\{p',p''\}$ satisfies
$$
p>\sigma,\quad p>{5\over2},\quad {5\over\varrho}-{5\over p}<1.
\leqno(3.42)
$$
Similarly we can prove that there exists $q$ such that
$$
q>\sigma,\quad q>5\quad {\rm and}\quad {5\over\varrho}-{5\over q}<2.
\leqno(3.43)
$$
Define $\bar p,\bar q$ such that ${1\over p}+{1\over\bar p}={1\over\sigma}$, 
${1\over q}+{1\over\bar q}={1\over\sigma}$.
Assume ${5\over3}<\sigma<\infty$. In view of Theorem 2.1 for a solution to 
problem (2.10) we have
$$\eqal{
&\|h\|_{W_\sigma^{2,1}(\Omega^t)}+\|\nabla q\|_{L_\sigma(\Omega^t)}\le c
(\|v\cdot\nabla h\|_{L_\sigma(\Omega^t)}\cr
&\quad+\|h\cdot\nabla v\|_{L_\sigma(\Omega^t)}+
\|\alpha_\theta\vartheta f\|_{L_\sigma(\Omega^t)}+
\|\alpha g\|_{L_\sigma(\Omega^t)}\cr
&\quad+\|h(0)\|_{W_\sigma^{2-2/\sigma}(\Omega)}).\cr}
\leqno(3.44)
$$
By (3.42) and (3.43) we estimate the first term on the r.h.s. of (3.44) by
$$\eqal{
&\|v\|_{L_q(\Omega^t)}\|\nabla h\|_{L_{\bar q}(\Omega^t)}\le c
\|v\|_{W_\varrho^{2,1}(\Omega^t)}(\varepsilon_2
\|h\|_{W_\sigma^{2,1}(\Omega^t)}\cr
&\quad+c(\varepsilon_2)\|h\|_{L_2(\Omega^t)})\cr}
\leqno(3.45)
$$
the second by
$$
\|\nabla v\|_{L_p(\Omega^t)}\|h\|_{L_{\bar p}(\Omega^t)}\le c
\|v\|_{W_\varrho^{2,1}(\Omega^t)}(\varepsilon_3)\|h\|_{W_\sigma^{2,1}(\Omega^t)}
+c(\varepsilon_3)\|h\|_{L_2(\Omega^t)})
\leqno(3.46)
$$
the third by
$$
c_1\|f\|_{L_\infty(\Omega^t0}(\varepsilon_4
\|\vartheta\|_{W_\sigma^{2,1}(\Omega^t)}+c(\varepsilon_4)
\|\vartheta\|_{L_2(\Omega^t)})
\leqno(3.47)
$$
the fourth by
$$
c_1\|g\|_{L_\sigma(\Omega^t)}.
\leqno(3.48)
$$
In view of [3, Ch. 4, Sect. 9, Th. 9.1] for any solution to problem (2.13) 
we have
$$\eqal{
&\|\vartheta\|_{W_\sigma^{2,1}(\Omega^t)}\le c
(\|v\cdot\nabla\vartheta\|_{L_\sigma(\Omega^t)}\cr
&\quad+\|h\cdot\nabla\theta\|_{L_\sigma(\Omega^t)}+
\|\vartheta(0)\|_{W_\sigma^{2-2/\sigma}(\Omega^t)}).\cr}
\leqno(3.49)
$$
By (3.42) and (3.43) we estimate the first term on the r.h.s. of (3.49) by
$$
c\|v\|_{W_\varrho^{2,1}(\Omega^t)}(\varepsilon_5
\|\vartheta\|_{W_\sigma^{2,1}(\Omega^t)}+c(\varepsilon_5)
\|\vartheta\|_{L_2(\Omega^t)})
\leqno(3.50)
$$
the second by
$$
c\|\theta\|_{W_\varrho^{2,1}(\Omega^t)}(\varepsilon_6
\|h\|_{W_\sigma^{2,1}(\Omega^t)}+c(\varepsilon_6)
\|h\|_{L_2(\Omega^t)}).
\leqno(3.51)
$$
We choose $r$ such that ${5\over3}<r<{10\over3}$ and $r\le\sigma$. By (3.1), 
the imbedding
$$
\|h\|_{L_\infty(0,t;L_3(\Omega))}\le c\|h\|_{W_r^{2,1}(\Omega^t)}
$$
and (3.24) there exists a constant $B_3$ depending on constants from 
imbedding theorems and data such that
$$
\|h\|_{L_2(\Omega^t)}+\|\vartheta\|_{L_2(\Omega^t)}\le B_3.
\leqno(3.52)
$$
Assuming that $\varepsilon_2-\varepsilon_6$ are sufficiently small and using 
(3.44)--(3.52) we obtain
$$
\|h\|_{W_\sigma^{2,1}(\Omega^t)}+\|\nabla q\|_{L_\sigma(\Omega^t)}+
\|\vartheta\|_{W_\sigma^{2,1}(\Omega^t)}\le B_4,
\leqno(3.53)
$$
where $B_4$ is some constant depending on data. By (3.40), (3.41) and (3.53) 
the proof is finished.
\kwadrat

\section{References}

\item{1.} Alame, W.: On existence of solutions for the nonstationary Stokes 
system with slip boundary conditions, Appl. Math. 32 (2) (2005), 195--223.

\item{2.} Besov, O. V.; Il'in, V. P.; Nikol'skii, S. M.: Integral 
representation of functions and imbedding theorems, Nauka, Moscow 1975 
(in Russian).

\item{3.} Ladyzhenskaya, O. A.; Solonnikov, V. A.; Uraltseva, N. N.: Linear 
and quasilinear equations of parabolic type, Moscow 1967 (in Russian).

\item{4.} Nowakowski, B.; Zaj\c aczkowski, W. M.: Very weak solutions to the 
boundary-value problem of the homogeneous hest equation in bounded domains 
(to be published).

\item{5.} Soca\l a, J.; Zaj\c aczkowski, W. M.: Long time existence of 
solutions to 2d Navier-Stokes equations with heat convection, Appl. Math.
36 (4) (2009), 453--463.

\item{6.} Soca\l a, J.; Zaj\c aczkowski, W. M.: Long time existence of 
regular solutions to 3d Navier-Stokes equations coupled with the heat 
convection (to be published).

\item{7.} Renc\l awowicz, J.; Zaj\c aczkowski, W. M.: Large time regular 
solutions to the Navier-\ -Stokes equations in cylindrical domains, TMNA 
32 (2008), 69--87.

\item{8.} Zaj\c aczkowski, W. M.: Long time existence of regular solutions 
to Navier-Stokes equations in cylindrical domains under boundary slip 
conditions, Studia Math. 169 (3) (2005), 243--285.

\item{9.} Zaj\c aczkowski, W. M.: Global special solutions to the 
Navier-Stokes equations in a cylindrical domain without the axis of 
symmetry, TMNA 24 (2004), 69--105.

\item{10.} Zaj\c aczkowski, W. M.: Global existence of axially symmetric 
solutions of incompressible Navier-Stokes equations with large angular 
component of velocity, Colloq. Math. 100 (2004), 243--263.

\bye